\documentclass{article}
\usepackage{amssymb,latexsym,amsmath,amsthm,mathrsfs}
\usepackage{graphicx}
\newcommand{\comment}[1]{}

\begin{document}
\title{On the partition of numbers into parts of a given type and number\footnote{Presented to the
St. Petersburg Academy on August 18, 1768.
Originally published as
{\em De partitione numerorum in partes tam numero quam specie datas},
Novi Commentarii academiae scientiarum Petropolitanae \textbf{14} (1770), 168--187.
E394 in the Enestr{\"o}m index.
Translated from the Latin by Jordan Bell,
Department of Mathematics, University of Toronto, Toronto, Canada.
Email: jordan.bell@utoronto.ca}}
\author{Leonhard Euler}
\date{}
\maketitle

1. Some time ago I had treated the problem of the partition of numbers,
 which looks for how many different ways a given number can be 
separated into two, three, or generally any number, of parts. I had been
careful that I present nothing about this solution by induction, as
we often use in solving this kind of problem. 
The method which I used seems such that it could be applied
with equal success to other problems, in particular
the common problem which searches for how many ways a given
number can be thrown as a given number of dice, which indeed I have
decided to explain here in a way that can be easily
generalized.

2. Since we are looking for how many ways a given number $N$ can occur
by throwing a given number $n$ of dice, here the question reduces to this,
how many distinct ways can
can a given number $N$ be resolved into $n$ parts,
which are altogether $1,2,3,4,5$ or $6$, if the sides of the dice are
marked with these numbers.
From this a more general question presents itself,
in how many distinct ways can a given number $N$ be divided into $n$ parts,
which are altogether $\alpha,\beta,\gamma,\delta$ etc.,  
the number of which is $=m$, such that both the number and type
of parts, in which a given number shall be resolved, are given.

3. Namely let the dice, not just in this particular case of six sides, but indeed
have $m$ sides or faces, so that in each the sides are
marked with the numbers
$\alpha,\beta,\gamma,\delta$, etc., and one then asks how many ways a given
number $N$
can be produced by throwing $n$ of these dice.
It could also be assumed that the dice are different from each other,
so that each would have a particular number of faces,
which would moreover be inscribed with particular
numbers;
truly the solution of this question, which I have generalized
from the common problem of dice, can found without excessive difficulty.

4. Indeed, I consider the numbers, which the sides of the dice are marked with,
as exponents of some quantity $x$, so that for common dice we would have
this expression
\[
x^1+x^2+x^3+x^4+x^5+x^6,
\]
where I take unity as the coefficient of each power,
since each number is designated by the exponent equal to it.
But if we take the square of this expression,
each power of $x$ will receive a coefficient that indicates how many ways
this power can result from the multiplication of two terms of
the expression, that is, how many ways its exponent can be produced
from the addition of two numbers from the sequence $1,2,3,4,5,6$.
Therefore by expanding the square of our expression, if the term
$Mx^N$ occurs in it, it follows from this that the number $N$ can be 
thrown as two dice
in as many ways as the coefficient
$M$ contains unities.

5. In a similar way, it is clear that if one takes
the cube of this expression $(x+x^2+x^3+x^4+x^5+x^6)^3$,
in its expansion any particular power $x^N$ occurs just as many
times as the ways in which the exponent $N$ can arise from the addition
of three numbers from the sequence $1,2,3,4,5,6$;
thus if $M$ is the coefficient of a power, and the whole term is $Mx^N$,
we conclude from this that the the number $N$ can be produced from
throws of three dice in as many ways as the coefficient $M$
contains unities. Therefore in general, if $n$ is taken as the exponent
of our expression $(x+x^2+x^3+x^4+x^5+x^6)^n$,
by expanding according to powers of $x$, any term $Mx^N$
shows that, if the number of dice $=n$,
the number $N$ can be made by throwing them together in as many
ways as the coefficient $M$ contains unities.

6. Therefore if the number of dice were $=n$ and we search for
how many ways a given number $N$ can be made by throwing these,
the question is resolved by the expansion of this formula
\[
(x+x^2+x^3+x^4+x^5+x^6)^n;
\]
now since the first term will be $x^n$, the last indeed $x^{6n}$,
the progression of the terms follows in this way
\[
x^n+Ax^{n+1}+Bx^{n+2}+Cx^{n+3}+\cdots+Mx^N+\cdots+x^{6n},
\]
where for any term $Mx^N$ it is clear that the number $N$ which equals
the exponent can be made in as many ways as the coefficient $M$ contains
unities;
from this it is immediately clear that the question does not even occur unless
the proposed number $N$ is contained within the limits $n$ and $6n$.
The whole matter thus now becomes that this progression, or the coefficients
of all its terms, be assigned.

7. Therefore to find these, let the formula that is to be expanded
be represented in this way 
\[
x^n(1+x+x^2+x^3+x^4+x^5)^n=V,
\]
then indeed for its expansion let 
\[
V=x^n(1+Ax+Bx^2+Cx^3+Dx^4+Ex^5+Fx^6+\textrm{etc.}).
\]
And by putting $\frac{V}{x^n}=Z$ by differentiating the logarithm 
it will be for the first
\[
\frac{xdZ}{Zdx}=\frac{nx+2nx^2+3nx^3+4nx^4+5nx^5}{1+x+x^2+x^3+x^4+x^5}.
\]
As well, the value for the latter will follow
\[
\frac{xdZ}{Zdx}=\frac{Ax+2Bx^2+3Cx^3+4Dx^4+5Ex^5+6Fx^6+\textrm{etc.}}{1+Ax+Bx^2+Cx^3+Dx^4+Ex^5+Fx^6+\textrm{etc.}},
\]
and since these two expressions are equal to each other,
from this the values of the coefficients can be determined.

8. Having constituted the equality of these two expressions,
this equation follows
\[
\begin{split}
&\begin{array}{lllllllll}
nx&+nAx^2&+nBx^3&+nCx^4&+nDx^5&+nEx^6&+nFx^7&+nGx^8&+\textrm{etc.}\\
&+2n&+2nA&+2nB&+2nC&+2nD&+2nE&+2nF&\\
&&+3n&+3nA&+3nB&+3nC&+3nD&+3nE&\\
&&&+4n&+4nA&+4nB&+4nC&+4nD&\\
&&&&+5n&+5nA&+5nB&+5nC&
\end{array}\\
&=
\begin{array}{lllllllll}
Ax&+2Bx^2&+3Cx^3&+4Dx^4&+5Ex^5&+6Fx^6&+7Gx^7&+8Hx^8&+\textrm{etc.}\\
&+A&+2B&+3C&+4D&+5E&+6F&+7G&\\
&&+A&+2B&+3C&+4D&+5E&+6F&\\
&&&+A&+2B&+3C&+4D&+5E&\\
&&&&+A&+2B&+3C&+4D&\\
&&&&&+A&+2B&+3C&
\end{array}
\end{split}
\]
and these two expressions, since all their terms are equal to each other,
will provide the values of all the coefficients.

9. The following determinations are now obtained:
\begin{eqnarray*}
A&=&n\\
2B&=&(n-1)A+2n,\\
3C&=&(n-2)B+(2n-1)A+3n,\\
4D&=&(n-3)C+(2n-2)B+(3n-1)A+4n,\\
5E&=&(n-4)D+(2n-3)C+(3n-2)+(4n-1)A+5n,\\
6F&=&(n-5)E+(2n-4)D+(3n-3)C+(4n-2)B+(5n-1)A,\\
7G&=&(n-6)F+(2n-5)E+(3n-4)D+(4n-3)C+(5n-2)B,\\
8H&=&(n-7)G+(2n-6)F+(3n-5)E+(4n-4)D+(5n-3)C\\
&&\textrm{etc.}
\end{eqnarray*}
Thus any coefficient can be determined by the preceding five, and with these found it will be
\[
V=x^n+Ax^{n+1}+Bx^{n+2}+Cx^{n+3}+Dx^{n+4}+Ex^{n+5}+\textrm{etc.}
\]
and this is the general solution to the problem for $n$ faces.

10. If each of the preceding is subtracted from the above equations,
the following much simpler determinations can be obtained:
\begin{eqnarray*}
A&=&n,\\
2B&=&nA+n,\\
3C&=&nB+nA+n,\\
4D&=&nC+nB+nA+n,\\
5E&=&nD+nC+nB+nA+n,\\
6F&=&nE+nD+nC+nB+nA-5n,\\
7G&=&nF+nE+nD+nC+nB-(5n-1)A,\\
8H&=&nG+nF+nE+nD+nC-(5n-2)B\\
&&\textrm{etc.}
\end{eqnarray*}
If again the differences are taken, these relations will become even simpler, in this way:
\begin{eqnarray*}
2B&=&(n+1)A,\\
3C&=&(n+2)B,\\
4D&=&(n+3)C,\\
5E&=&(n+4)D,\\
6F&=&(n+5)E-6n,\\
7G&=&(n+6)F-(6n-1)A+5n,\\
8H&=&(n+7)G-(6n-2)B+(5n-1)A,\\
9I&=&(n+8)H-(6n-3)C+(5n-2)B,\\
10K&=&(n+9)I-(6n-4)D+(5n-3)C\\
&&\textrm{etc.}
\end{eqnarray*}

11. Then, if the number of dice were $2,3$ or $4$,
the law for the progression of the coefficients will be as follows,
\[
\begin{array}{l|l|l}
\textrm{for two}&\textrm{for three}&\textrm{for four}\\
A=2&3&4\\
2B=3A&4A&5A\\
3C=4B&5B&6B\\
4D=5C&6C&7C\\
5E=6D&7D&8D\\
6F=7E-12&8E-18&9E-24\\
7G=8F-11A+10&9F-17A+15&10F-23A+20\\
8H=9G-10B+9A&10G-16B+14A&11G-22B+19A\\
9I=10H-9C+8B&11H-15C+13B&12H-21C+18B\\
10K=11I-8D+7C&12I-14D+12C&13I-20D+17C\\
11L=12K-7E+6D&13K-13E+11D&14K-19E+16D\\
12M=13L-6F+5E&14L-12F+10E&15L-18F+15E\\
\textrm{etc.}&\textrm{etc.}&\textrm{etc.}
\end{array}
\]
Therefore any coefficient can be determined by three preceding, where it
is particularly noteworthy here that finally they shall abate to nothing, and 
the same happens for those following the first to disappear, which
is not entirely clear from this law.

12. To help us clearly understand this law, let the formula
\[
(N)^{(n)}
\]  
denote the number of cases in which the number $N$ can be made
by $n$ dice, so that it would be
\[
\begin{split}
&(n)^{(n)}=1,\quad (n+1)^{(n)}=A,\quad (n+2)^{(n)}=B,\quad (n+3)^{(n)}=C,\\
&(n+4)^{(n)}=D,\ldots,(n+9)^{(n)}=I \quad \textrm{and} \quad (n+10)^{(n)}=K.
\end{split}
\]
It will therefore be
\[
10(n+10)^{(n)}=(n+9)(n+9)^{(n)}-(6n-4)(n+4)^{(n)}+(5n-3)(n+3)^{(n)},
\]
and here it can be concluded generally to be
\[
\begin{split}
&\lambda(n+\lambda)^{(n)}=(n+\lambda-1)(n+\lambda-1)^{(n)}-(6n+6-\lambda)(n+\lambda-6)^{(n)}\\
&+(5n+7-\lambda)(n+\lambda-7)^{(n)}.
\end{split}
\]
Let us now put $n+\lambda=N$, and so that $\lambda=N-n$, and it will be
\[
(N)^{(n)}=\frac{(N-1)(N-1)^{(n)}-(7n+6-N)(N-6)^{(n)}+(6n+7-N)(N-7)^{(n)}}{N-n},
\]
where it is noted to always be $(P)^{(n)}=0$ when $P<n$.

13. Also, these coefficients can easily be defined for any number of dice,
if those for the number of dice one less are known. For if it were
\[
\begin{split}
&(x+x^2+x^3+x^4+x^5+x^6)^n\\
&=x^n+Ax^{n+1}+Bx^{n+2}+Cx^{n+3}+Dx^{n+4}+\textrm{etc.}
\end{split}
\]
let it be put
\[
\begin{split}
&(x+x^2+x^3+x^4+x^5+x^6)^{n+1}\\
&x^{n+1}+A'x^{n+2}+B'x^{n+3}+C'x^{n+4}+D'x^{n+5}+\textrm{etc.}
\end{split}
\]
then it will be, because this expression is equal to the former
multiplied by $x+x^2+x^3+x^4+x^5+x^6$,
\[
\begin{array}{l|l}
A'=A+1&\textrm{then taking the differentials}\\
B'=B+A+1&B'=A'+B\\
C'=C+B+A+1&C'=B'+C\\
D'=D+C+B+A+1&D'=C'+D\\
E'=E+D+C+B+A+1&E'=D'+E\\
F'=F+E+D+C+B+A&F'=E'+F-1\\
G'=G+F+E+D+C+B&G'=F'+G-A\\
\textrm{etc.}&\textrm{etc.}
\end{array}
\]

14. Thus if the way of denoting introduced above is used,
from the equation $G'=F'+G-A$ it is
\[
(n+8)^{(n+1)}=(n+7)^{(n+1)}+(n+7)^{(n)}-(n+1)^{(n)},
\]
which will be represented in general as
\[
(n+1+\lambda)^{(n+1)}=(n+\lambda)^{(n+1)}+(n+\lambda)^{(n)}-(n+\lambda-6)^{(n)}.
\]
Now if $n+\lambda$ is written for $N$, it will be
\[
(N+1)^{(n+1)}=(N)^{(n+1)}+(N)^{(n)}-(N-6)^{(n)},
\]
where it should be noted that whenever $N-6<n$ it will be $(N-6)^{(n)}=0$.
It is clear immediately that all these numbers are integers, which was less apparent
from the previous law.

\begin{table}
\caption{Exhibiting how many ways a number $N$ can be made
by $n$ dice}
\[
\begin{array}{r|r|r|r|r|r|r|r|r}
N&n=1&n=2&n=3&n=4&n=5&n=6&n=7&n=8\\
1&1&0&0&0&0&0&0&0\\
2&1&1&0&0&0&0&0&0\\
3&1&2&1&0&0&0&0&0\\
4&1&3&3&1&0&0&0&0\\
5&1&4&6&4&1&0&0&0\\
6&1&5&10&10&5&1&0&0\\
7&0&6&15&20&15&6&1&0\\
8&0&5&21&35&35&21&7&1\\
9&0&4&25&56&70&56&28&8\\
10&0&3&27&80&126&126&84&36\\
11&0&2&27&104&205&252&210&120\\
12&0&1&25&125&305&456&462&330\\
13&0&0&21&140&420&756&917&792\\
14&0&0&15&146&540&1161&1667&1708\\
15&0&0&10&140&651&1666&2807&3368\\
16&0&0&6&125&735&2247&4417&6147\\
17&0&0&3&104&780&2856&6538&10480\\
18&0&0&1&80&780&3431&9142&16808\\
19&0&0&0&56&735&3906&12117&25488\\
20&0&0&0&35&651&4221&15267&36688\\
21&0&0&0&20&540&4332&18327&50288\\
22&0&0&0&10&420&4221&20993&65808\\
23&0&0&0&4&305&3906&22967&82384\\
24&0&0&0&1&205&3431&24017&98813\\
25&0&0&0&0&126&2856&24017&113688\\
26&0&0&0&0&70&2247&22967&12588\\
27&0&0&0&0&35&1666&20993&133288\\
28&0&0&0&0&15&1161&18327&135954\\
29&0&0&0&0&5&756&15267&133288\\
30&0&0&0&0&1&456&12117&125588\\
31&0&0&0&0&0&252&9142&113688\\
32&0&0&0&0&0&126&6538&98813\\
33&0&0&0&0&0&56&4417&82384\\
34&0&0&0&0&0&21&2807&65808\\
35&0&0&0&0&0&6&1667&50288\\
36&0&0&0&0&0&1&917&36688
\end{array}
\]
\end{table}

15. In these series also the properties found in \S 12 occur;
thus if $n=6$, it will be
\[
(N)^{(6)}=\frac{(N-1)(N-1)^{(6)}-(48-N)(N-6)^{(6)}+(43-N)(N-7)^{(6)}}{N-6},
\]
where, if for example $N=25$, it will be
\[
(25)^{(6)}=\frac{24\cdot (24)^{(6)}-23\cdot (19)^{(6)}+18\cdot (18)^{(6)}}{19};
\]
and it is $(24)^{(6)}=3431,(19)^{(6)}=3906,(18)^{(6)}=3431$, 
hence
\[
(25)^{(6)}=\frac{24\cdot 3431-23\cdot 3906+18\cdot 3431}{19}=\frac{54264}{19}=2856,
\]
as is obtained in the table. Similarly if $N=29$, it will be
\[
(29)^{(6)}=\frac{28\cdot (28)^{(6)}-19\cdot (23)^{(6)}+14\cdot (22)^{(6)}}{23};
\]
then because $(28)^{(6)}=1161,(23)^{(6)}=3906$ and $(22)^{(6)}=4221$ it will
be
\[
(29)^{(6)}=\frac{32508-74214+59094}{23}=\frac{17388}{23}=756.
\]

16. Truly the expansion of the formula $V$ (\S 7) can be done in a different
way, such that any term can be completely assigned without
the use of the preceding. For since it is
\[
1+x+x^2+x^3+x^4+x^5=\frac{1-x^6}{1-x},
\]
it will be
\[
V=\frac{x^n(1-x^6)^n}{(1-x)^n}
\]
and by expanding 
\begin{eqnarray*}
(1-x^6)^n&=&1-\frac{n}{1}x^6+\frac{n(n-1)}{1\cdot 2}x^{12}-\frac{n(n-1)(n-2)}{1\cdot 2\cdot 3}x^{18}\\
&&+\frac{n(n-1)(n-2)(n-3)}{1\cdot 2\cdot 3\cdot 4}x^{24}-\textrm{etc.},
\end{eqnarray*}
\begin{eqnarray*}
\frac{x^n}{(1-x)^n}&=&x^n+\frac{n}{1}x^{n+1}+\frac{n(n+1)}{1\cdot 2}x^{n+2}+\frac{n(n+1)(n+2)}{1\cdot 2\cdot 3}x^{n+3}\\
&&+\frac{n(n+1)(n+2)(n+3)}{1\cdot 2\cdot 3\cdot 4}x^{n+4}+\textrm{etc.}
\end{eqnarray*}
and thus it can be concluded to be
\begin{eqnarray*}
(n)^{(n)}&=&1,\\
(n+1)^{(n)}&=&\frac{n}{1},\\
(n+2)^{(n)}&=&\frac{n(n+1)}{1\cdot 2},\\
(n+3)^{(n)}&=&\frac{n(n+1)(n+2)}{1\cdot 2\cdot 3},\\
(n+4)^{(n)}&=&\frac{n(n+1)\cdots (n+3)}{1\cdot 2\cdot 3\cdot 4},\\
(n+5)^{(n)}&=&\frac{n(n+1)\cdots (n+4)}{1\cdot 2\cdot 3\cdot 4\cdot 5},\\
(n+6)^{(n)}&=&\frac{n(n+1)\cdots (n+5)}{1\cdot 2\cdot 3\cdots 6}-\frac{n}{1}\cdot 1,\\
(n+7)^{(n)}&=&\frac{n(n+1)\cdots (n+6)}{1\cdot 2\cdot 3\cdots 7}-\frac{n}{1}\cdot \frac{n}{1},\\
(n+8)^{(n)}&=&\frac{n(n+1)\cdots (n+7)}{1\cdot 2\cdot 3\cdots 8}-\frac{n}{1}\cdot \frac{n(n+1)}{1\cdot 2},\\
(n+9)^{(n)}&=&\frac{n(n+1)\cdots (n+8)}{1\cdot 2\cdot 3\cdots 9}-\frac{n}{1}\cdot \frac{n(n+1)(n+2)}{1\cdot 2\cdot 3},\\
(n+10)^{(n)}&=&\frac{n(n+1)\cdots (n+9)}{1\cdot 2\cdot 3\cdots 10}-\frac{n}{1}\cdot \frac{n(n+1)\cdots (n+3)}{1\cdot 2\cdot 3\cdot 4},\\
(n+11)^{(n)}&=&\frac{n(n+1)\cdots (n+10)}{1\cdot 2\cdot 3\cdots 11}-\frac{n}{1}\cdot \frac{n(n+1)\cdots (n+4)}{1\cdot 2\cdot 3\cdot 4\cdot 5},\\
(n+12)^{(n)}&=&\frac{n(n+1)\cdots (n+11)}{1\cdot 2\cdot 3\cdots 12}-\frac{n}{1}\cdot \frac{n(n+1)\cdots (n+5)}{1\cdot 2\cdot 3\cdots 6}+\frac{n(n-1)}{1\cdot 2}\cdot 1,\\
(n+13)^{(n)}&=&\frac{n(n+1)\cdots (n+12)}{1\cdot 2\cdot 3\cdots 13}-\frac{n}{1}\cdot \frac{n(n+1)\cdots (n+6)}{1\cdot 2\cdot 3\cdots 7}+\frac{n(n-1)}{1\cdot 2}\cdot \frac{n}{1}\\
&&\textrm{etc.,}
\end{eqnarray*}
from which it can concluded in general
\[
\begin{split}
&(n+\lambda)^{(n)}\\
&=\frac{n(n+1)\cdots (n+\lambda-1)}{1\cdot 2\cdot 3\cdots \lambda}-
\frac{n}{1}\cdot \frac{n(n+1)\cdots (n+\lambda-7)}{1\cdot 2\cdot 3\cdots (\lambda-6)}+
\frac{n(n-1)}{1\cdot 2}\cdot \frac{n(n+1)\cdots (n+\lambda-13)}{1\cdot 2\cdot 3\cdots (\lambda-12)}\\
&-\frac{n(n-1)(n-2)}{1\cdot 2\cdot 3}\cdot \frac{n(n+1)\cdots (n+\lambda-19)}{1\cdot 2\cdot 3\cdots (\lambda-18)}+\frac{n(n-1)(n-2)(n-3)}{1\cdot 2\cdot 3\cdot 4}\cdot
\frac{n(n+1)\cdots (n+\lambda-25)}{1\cdot 2\cdot 3\cdots (\lambda-24)}\\
&-\textrm{etc.}
\end{split}
\]

17. This solution can be adapted to dice possessing any other number of faces.
For let the number of faces on all the dice be $m$, which shall
be marked with the numbers $1,2,3,\ldots, m$, 
while the number of dice itself shall be $=n$ for
which the number of throws are sought in which the given number $N$ land.
Or, which ends up being the same thing, it is sought how many ways the number $N$
can be resolved into $n$ parts, each of which is comprised from the sequence
of
numbers $1,2,3,\ldots,m$;
here indeed it should be noted that not only different partitions but also
different orders of the same parts are counted, as usually happens in dice,
where for example the throws $3,4$ and $4,3$ are considered to be two different
cases.

18. But if therefore this symbol $(N)^{(n)}$ denotes the number of cases
in which the number $N$ can be produced by throwing $n$ dice,
each of which have $m$ sides marked with the numbers $1,2,3,\ldots,m$,
it should first be noted that $(n)^{(n)}=1$, and if $N<n$ then $(N)^{(n)}=0$.
Next, if $N=mn$ then too $(mn)^{(n)}=1$, and
if $N>mn$ it will be $(N)^{(n)}=0$.
Finally, if either $N=n+\lambda$ or $N=mn-\lambda$, the number of
cases is the same, namely
\[
(n+\lambda)^{(n)}=(mn-\lambda)^{(n)}.
\]
The first formula can also be developed as
\begin{eqnarray*}
(n+\lambda)^{(n)}&=&\frac{n(n+1)\cdots(n+\lambda-1)}{1\cdot 2\cdot 3\cdots \lambda}-\frac{n}{1}\cdot \frac{n(n+1)\cdots (n+\lambda-m-1)}{1\cdot 2\cdot 3\cdots (\lambda-m)}\\
&&+\frac{n(n-1)}{1\cdot 2}\cdot \frac{n(n+1)\cdots (n+\lambda-2m-1)}{1\cdot 2\cdot 3\cdots (\lambda-2m)}\\
&&-\frac{n(n-1)(n-2)}{1\cdot 2\cdot 3}\cdot
\frac{n(n+1)\cdots (n+\lambda-3m-1)}{1\cdot 2\cdot 3\cdots (\lambda-3m)}+\textrm{etc.}
\end{eqnarray*}

19. As well, these numbers can very easily be determined from the preceding,
where the number of dice is one less. For in general it will be,
if the number of all the dice were $=m$ and each were
marked with the numbers $1,2,3,\ldots,m$,
\[
(N+1)^{(n+1)}=(N)^{(n+1)}+(N)^{(n)}-(N-m)^{(n)}
\]
or
\[
(N+1)^{(n)}=(N)^{(n)}+(N)^{(n-1)}-(N-m)^{(n-1)}.
\]
Then, if $n+\lambda$ is written for $N+1$, we will obtain
\[
(n+\lambda)^{(n)}=(n+\lambda-1)^{(n)}+(n+\lambda-1)^{(n-1)}-(n+\lambda-m-1)^{(n-1)}.
\]
Finally, for the same number of dice, this number $n$ depends
thus on preceding numbers
\begin{eqnarray*}
\lambda(n+\lambda)^{(n)}&=&(n+\lambda-1)(n+\lambda-1)^{(n)}-(mn+m-\lambda)(n+\lambda-m)^{(n)}\\
&&+(mn-n+m+1-\lambda)(n+\lambda-m-1)^{(n)}.
\end{eqnarray*}
Again it is noted that the sum of all these numbers is $=m^n$.

20. The question can be resolved in a similar way even if not all the 
dice occur with the same number of faces. Let us take three dice,
the first a hexahedron bearing the numbers
$1,2,3,4,5,6$, the second an octahedron bearing the numbers $1,2,3,\ldots,8$,
and the third a dodecahedron bearing the numbers $1,2,3,\ldots,12$;
if we now ask for how many ways a given number $N$ can be decomposed,
we should expand this product
\[
(x+x^2+x^3+\cdots+x^6)(x+x^2+x^3+\cdots+x^8)(x+x^2+x^3+\cdots+x^{12})=V,
\]
and the coefficient of the power $x^N$ will display the number of cases.
Now, since here it will be
\[
V=\frac{x^3(1-x^6)(1-x^8)(1-x^{12})}{(1-x)^3},
\]
by expanding the numerator it will be
\[
V=\frac{x^3-x^9-x^{11}-x^{15}+x^{17}+x^{21}+x^{23}-x^{29}}{(1-x)^3}.
\]

21. Here the numerator is multiplied by $\frac{1}{(1-x)^3}$,
or the series
\[
1+3x+6x^2+10x^3+15x^4+21x^5+28x^6+36x^7+\textrm{etc.},
\]
whose coefficients are the triangular numbers; 
and since the $n$th triangular number is $\frac{n(n+1)}{2}$,
any term of this series will be
\[
\frac{n(n+1)}{2}x^{n-1} \quad \textrm{or} \quad \frac{(n-1)(n-2)}{2}x^{n-3}.
\]
Now multiplying by the numerator, the coefficient of
the power $x^n$ turns out to be
\[
\begin{split}
&\frac{(n-1)(n-2)}{2}-\frac{(n-7)(n-8)}{2}-\frac{(n-9)(n-10)}{2}-\frac{(n-13)(n-14)}{2}\\
&+\frac{(n-15)(n-16)}{2}+\frac{(n-19)(n-20)}{2}+\frac{(n-21)(n-22)}{2}+\frac{(n-27)(n-28)}{2}.
\end{split}
\]
This expression does not need to be continued beyond the case when negative
factors are arrived at. 

22. Noticing however for the denominator $(1-x)^3=1-3x+3x^2-x^3$,
the series that is sought will be recurrent, arising from the ladder relation
$3,-3,+1$,
providing that we have the rule for the terms of the numerator.
Then the following coefficients for each exponent are found:
\[
\begin{array}{c|c|c|c}
\textrm{Exponents}&\textrm{Coefficients}&\textrm{Exponents}&\textrm{Coefficients}\\
\hline
3&1&15&47\\
4&3&16&45\\
5&6&17&42\\
6&10&18&38\\
7&15&19&33\\
8&21&20&27\\
9&27&21&21\\
10&33&22&15\\
11&38&23&10\\
12&42&24&6\\
13&45&25&3\\
14&47&26&1
\end{array}
\]  
The numbers here are not greater than $26$, since $26=6+8+12$,
and the sum of all the cases is $576=6\cdot 8\cdot 12$.

23. Since here we have succeeded in the resolution of numbers into parts
of a given number and type without the help of induction, some elegant
Theorems of Fermat come to my mind; while they have not yet been demonstrated,
it seems that this method will perhaps lead to demonstrations of
them. 
For Fermat had asserted that all numbers are either triangles or
the aggregate of two or three triangles, and because zero also
occurs in the order of triangles, the theorem can thus be stated 
as that all numbers can be resolved into three triangles.
Thus if the triangular numbers are taken for exponents and this series
is formed
\[
1+x^1+x^3+x^6+x^{10}+x^{15}+x^{21}+x^{28}+\textrm{etc.}=S,
\]
it needs to be demonstrated that if the cube of this series is expanded,
then every power of $x$ occurs without omission;
for if this could be demonstrated, the demonstration of this Theorem
of Fermat would be possessed.

24. In a similar way, if we take the fourth power of this series
\[
1+x^1+x^4+x^9+x^{16}+x^{25}+x^{36}+\textrm{etc.}=S,
\]
and we could show that every power of $x$ occurs in it, then we would
have a demonstration of the Theorem of Fermat, that all numbers are made
from the addition of four squares.

Indeed in general if we put
\[
S=1+x^1+x^m+x^{3m-3}+x^{6m-8}+x^{10m-15}+x^{15m-24}+x^{21m-35}+\textrm{etc.}
\]
and we take the power of this series of the exponent $m$, it needs
to be demonstrated that every power of $x$ will be produced in it,
and thus that all numbers are the aggregate of $m$ or less polygonal numbers
with a number of sides $=m$.

25. From these same principles another way presents itself for investigating
the demonstrations, which differs from the preceding in that where
there not only different parts but also different orders were considered,
here the stipulation about the order is omitted. 
Namely for the resolution into triangular numbers let this formula be
constituted
\[
\frac{1}{(1-z)(1-xz)(1-x^3z)(1-x^6z)(1-x^{10}z)(1-x^{15}z)\, \textrm{etc.}},
\] 
which expanded produces this series
\[
1+Pz+Qz^2+Rz^3+Sz^4+Tz^5+\textrm{etc.},
\]
such that $P,Q,R,S$ etc. are some functions of $x$. Now it will clearly
be
\[
P=1+x+x^3+x^6+x^{10}+x^{15}+x^{21}+\textrm{etc.};
\]
while $Q$ will also contain those powers of $x$ whose exponents are
the aggregate of two triangles. It therefore  needs to be be demonstrated
that every power of $x$ occurs in the function $R$.

26. In a similar way for the resolution of numbers into four squares,
let this fraction be expanded
\[
\frac{1}{(1-z)(1-xz)(1-x^4z)(1-x^9z)(1-x^{16}z)(1-x^{25}z)\, \textrm{etc.}};
\] 
then if it is changed into the form
\[
1+Pz+Qz^2+Rz^3+Sz^4+\textrm{etc.},
\]
it needs to be demonstrated that the function $S$ contains every power of
$x$. For instance, $P$ is equal to the series $1+x+x^4+x^9+x^{16}+\textrm{etc.}$
and $Q$ also contains those powers of $x$ whose exponents are the aggregate
of two squares,
in which series therefore many powers are still absent. Further,
in $R$ there are the additional powers whose exponents are the aggregate of
three squares, and then in $S$ those will occur whose exponents are the sum
of four, so that in $S$ all numbers should occur as exponents.

27. From this principle, we can define how many solutions those problems
admit which are commonly referred to by Arithmeticians as the {\em Regula Virginum}\footnote{Translator: In English, {\em Rule of the Virgins}. See
Volume II of Leonard Eugene Dickson's
{\em History of the Theory of Numbers}, under the indices
``Coeci rule'' and ``Virgins''.}. In this way, the problems are reduced here to finding numbers
$p,q,r,s,t$ etc. that satisfy these two conditions
\[
ap+bq+cr+ds+\textrm{etc.}=n,
\]
and
\[
\alpha p+\beta q+\gamma r+\delta s+\textrm{etc.}=\nu,
\]
and now the question is, how many solutions occur in positive integers;
where indeed the numbers $a,b,c,d$ etc., $n$ and $\alpha,\beta,\gamma,\delta$ etc., $\nu$ may be taken to be integers, as if this were not so, it can be easily
reduced to this. It is at once 
apparent that if 
two numbers $p$ and $q$ are given to be investigated, more than one solution
cannot be given. Commonly the problem
only occurs with positive integral numbers for $p$ and $q$.

28. Now, for defining the number of all the solutions in any particular
case, so that nothing is done by induction or by trials, let us consider
this expression
\[
\frac{1}{(1-x^a y^\alpha)(1-x^b y^\beta)(1-x^c y^\gamma)(1-x^d y^\delta)\, \textrm{etc.}}
\]
and let it be expanded, which yields such a series
\[
1+Ax^{\ldots}y^{\ldots}+Bx^{\ldots}y^{\ldots}+Cx^{\ldots}y^{\ldots}
\, \textrm{etc.};
\]
if the term $Nx^n y^\nu$ occurs in this, the coefficient $N$ will indicate 
the number of solutions; and if it turns out that this term does not
occur, this indicates that there will be no solution. Thus
the whole problem turns here on investigating the coefficient of the
term $x^n y^\nu$.

\end{document}